\documentclass{amsart}
\usepackage{latexsym}
\usepackage{amssymb}
\usepackage{amsmath}
\usepackage{color}

\newtheorem{thm}{Theorem}[section]
\newtheorem{lem}[thm]{Lemma}

\newtheorem{openquestion}[thm]{Open Question}
\newtheorem{prop}[thm]{Proposition}
\newtheorem{cor}[thm]{Corollary}

\theoremstyle{remark}
    \newtheorem*{rem}{Remark}
\theoremstyle{definition}
    
\newtheorem{ex}[thm]{Example}

\numberwithin{equation}{section}

\def\Z {{\mathbb Z}}

\def\T {{\tilde T}}
\def\B {{\mathcal B}}
\def\C {{\mathcal C}}

\def\H {{\mathcal H}}
\def\M {{\mathcal M}}

\def\a {{\mathfrak a}}
\def\l {\left\langle}
\def\r {\right\rangle}

\def\tr{\mathrm{tr}}
\def\red#1 {\textcolor{red}{#1 }}
\def\blue#1 {\textcolor{blue}{#1 }}

\title[a new integral basis]{A new integral basis for the centre of the Hecke algebra of type $A$}
\author{Andrew Francis }
\address{School of Computing and Mathematics, University of Western Sydney, NSW 1797, Australia}
\email[Andrew~Francis]{a.francis@uws.edu.au}
\author{Lenny Jones}
\address{Department of Mathematics, Shippensburg University, Pennsylvania, USA}
\email[Lenny~Jones]{lkjone@ship.edu}

\begin{document}
\maketitle

\begin{abstract}
We describe a recursive algorithm that produces an integral basis for the centre of the Hecke algebra of type $A$ consisting of linear combinations of monomial symmetric polynomials of Jucys--Murphy elements. We also discuss the existence of integral bases for the centre of the Hecke algebra that consist solely of monomial symmetric polynomials of Jucys--Murphy elements. Finally, for $n=3$, we show that only one such basis exists for the centre of the Hecke algebra, by proving that there are exactly four bases for the centre of the corresponding symmetric group algebra which consist solely of monomial symmetric polynomials of Jucys--Murphy elements.
\end{abstract}

\section{Introduction}

The centre of the Iwahori--Hecke algebra of type $A_{n-1}$ shares a number of characteristics of the centre of the corresponding group algebra $\Z S_n$.  For instance, it has the same dimension (the number of conjugacy classes in $S_n$), it has an integral basis specialising to class sums~\cite{GR97}, and like the centre of the group algebra, it has long been known to be equal to the set of symmetric polynomials in Jucys--Murphy elements when considered over the rational function field~\cite{DJ87,Mur83}.  This last result has been extended to the integral Iwahori--Hecke algebra~\cite{FG:DJconj2006}, mirroring G.~Murphy's result for the centre of the integral group algebra precisely.  A natural question is then whether the basis given by Murphy in~\cite{Mur83} generalizes to the Iwahori--Hecke algebra. Murphy's basis consists of a particularly nice set of monomial symmetric polynomials of Jucys--Murphy elements.  Unfortunately, the direct generalization of this result \cite[Theorem 3.6]{Mat99}, fails in the Iwahori--Hecke algebras of $S_n$ for at least $n=3$, $4$ and 5 (see Section \ref{sec:basis}).

In this paper, we discuss the possibility that some sort of
variation of Murphy's monomial symmetric polynomial basis might exist
for the Hecke algebra, and we provide evidence to indicate that integral bases of the centre of the Hecke algebra consisting solely of monomial symmetric polynomials in
Jucys--Murphy elements are rare (Section \ref{sec:monbases}).
In particular, in Section~\ref{sec:S3.monbasis} we show that there is \emph{exactly one} basis of monomial symmetric polynomials in Jucys--Murphy elements for the centre of the Hecke algebra of $S_3$, by proving that there are exactly four such bases for the centre of the corresponding group algebra.
Our main theorem, however, shows that there always exists an integral basis for the centre of the Hecke algebra consisting of integral linear combinations of monomial symmetric polynomials in Jucys--Murphy elements (Theorem \ref{thm:monbasis}).  This is the first such basis for the centre of the Hecke algebra given in terms of symmetric polynomials, and the second description of an integral basis after~\cite{GR97,Fmb}.

We begin with Definitions and Preliminary results (Sections~\ref{sec:defs} and~\ref{sec:prelims}).

\section{Definitions}\label{sec:defs}

\subsection{Partitions and compositions}\label{Comp}

A \emph{composition} $\lambda$ is a finite ordered set of positive integers.  If $\lambda=(\lambda_1,\dots,\lambda_r)$, the $\lambda_i$ are called the \emph{components} of $\lambda$, and the number of components in $\lambda$, denoted $\ell(\lambda)$, is called its \emph{length}.  If $\lambda$ is a composition we write $|\lambda|=\sum_{i=1}^r\lambda_i$.  If $|\lambda|=n$ we say $\lambda$ is a \emph{composition of $n$}, and we write $\lambda\vDash n$. We denote the \emph{empty composition} (the unique composition of zero) as $\emptyset$, and we define $|\emptyset|:=0$. A \emph{partition} of $n$ is a composition whose components are weakly decreasing from left to right. If $\lambda$ is a partition of $n$ we write $\lambda\vdash n$. Repeated components of a composition will sometimes be denoted by exponents; for instance $(1^3,2)$ is shorthand for $(1,1,1,2)$.
If $\lambda=(\lambda_1,\dots,\lambda_r)\vDash n$ then we define $\lambda-1$ to be the composition $(\lambda_1-1,\dots,\lambda_r-1)$ of $n-\ell(\lambda)$, excluding zero components. For example, if $\lambda=(3,1,4)\vDash 8$, then $\lambda-1=(2,3)\vDash 5$. For a given $n$, if a composition $\lambda$ satisfies $\ell(\lambda)+|\lambda|\le n$, we define $\overline{\lambda}$ to be the partition of $n$ obtained from the parts of $\lambda$ by adding 1 to all parts, and $n-\ell(\lambda)-|\lambda|$ copies of 1.
 Note that the construction of $\overline\lambda$ is relative to $n$: for $\lambda=(3,1,4)$, we have $\overline\lambda=(5,4,2,1)$, if $n=12$; while $\overline\lambda=(5,4,2,1,1,1)$, if $n=14$. In particular, if $\lambda\vdash n$, then $\overline{\lambda-1}=\lambda$.

If $\lambda=(\lambda_1,\lambda_2,\ldots ,\lambda_t)$ is any nonempty composition, then we define \[\lambda^{\prime}:=\left\{ \begin{array}{cl}(\lambda_1,\lambda_2,\ldots ,\lambda_{t-1})& \mbox{ if $\ell(\lambda)\ge 2$} \\
\emptyset & \mbox{ if $\ell(\lambda)=1$}.\\ \end{array}\right.\]
\indent Let $\Lambda$ be the set of all compositions of all positive integers together with $\emptyset$. For the sake of convenience in the sequel, we define recursively the following ordering on $\Lambda$.
 For any two compositions $\lambda$ and $\mu$ in $\Lambda$, we have

\begin{equation}\label{eq:def:compn.order}
 \lambda<\mu \Longleftrightarrow \left\{ \begin{array}{l}
 |\lambda|<|\mu|, \mbox{ or}\\
 0<|\lambda|=|\mu| \mbox{ and } \lambda^{\prime}<\mu^{\prime}.\\
 \end{array}\right.
\end{equation}

\subsection{The symmetric group}\label{sec:Sn}

The symmetric group $S_n$ is the group of permutations of $\{1,\dots,n\}$.  We will use the Coxeter presentation for $S_n$ given by the generators $S=\{s_i=(i\ i+1)\mid 1\le i\le n-1\}$, and subject to the relations {$s_i^2=1$ for each $i$}; $s_is_j=s_js_i$ if $|i-j|\ge 2$; and $s_is_{i+1}s_i=s_{i+1}s_is_{i+1}$ for $1\le i\le n-2$. An expression $w=s_{i_1}\dots s_{i_k}$ for $w\in S_n$ is said to be \emph{reduced} if $k$ is minimal.  In this case, we say the length of $w$, denoted $\ell(w)$, is $k$.

An element $w$ of $S_n$ is called \emph{increasing} if it can be written in the form $s_{i_1}s_{i_2}\cdots s_{i_k}$ with $1\le i_1<i_2<\cdots <i_k<n$. For such an element $w$, we say \emph{$w$ is of shape $\lambda$}, if $\lambda$ is the composition  of $k$ defined as follows. Without any rearrangement of the $s_{i_j}$, we rewrite $w=s_{i_1}s_{i_2}\cdots s_{i_k}$ as $\alpha_1\alpha_2\cdots \alpha_t$, where the $\alpha_i$ are disjoint cycles of maximal length. Then $\lambda=(\ell(\alpha_1),\ell(\alpha_2),\ldots,\ell(\alpha_t))$.
For example, $s_1s_3s_4s_6$ has shape $(1,2,1)$.
Note that there exists an increasing element $w$ in $S_n$ of shape $\lambda$ if and only if $|\lambda|+\ell(\lambda)=k+t\le n$.

The conjugacy classes of $S_n$ are indexed by partitions of $n$.  Denote the conjugacy class of cycle type $\lambda\vdash n$ by $C_\lambda$.  Note that $C_\lambda$ is the unique conjugacy class that contains a minimal element of shape $\lambda-1$.
From the alternative perspective, if $w$ is an increasing element of shape $\lambda$ in $S_n$, then $\ell(\lambda)+|\lambda|\le n$ and $w$ is a minimal element in the conjugacy class $C_{\overline\lambda}$.

\subsection{Symmetric polynomials}
A polynomial in the variables
$\{X_1,\dots, X_n\}$ is \emph{symmetric} if it is fixed by the action of $S_n$ on the indices of the variables.
Each partition $\lambda=(\lambda_1,\dots,\lambda_r)$ with $r\le n$ determines a \emph{monomial symmetric polynomial},
\[m_\lambda=m_\lambda(X_1,\dots,X_n):=\sum_{\sigma\in S_n} X_{\sigma(1)}^{\lambda_1}\dots X_{\sigma(r)}^{\lambda_r},\]
where $m_{\emptyset}:=1.$ Note that we also write $m_{\lambda_1, \lambda_2, \ldots, \lambda_r}$ for $m_{\lambda}$.

For example, \[{m_{2,1}} (X_1,X_2,X_3):=X_1^2X_2+X_1^2X_3+X_2^2X_1+X_2^2X_3+X_3^2X_1+X_3^2X_2.\]

Given a composition $\lambda=(\lambda_1, \lambda_2, \ldots, \lambda_r)$ with $r\le n$, we define the \emph{monomial quasi-symmetric} polynomial $p^{\lambda}$ as
\[p^{\lambda}(X_1,X_2,\ldots,X_n):=\sum X_{j_1}^{\lambda_1}X_{j_2}^{\lambda_2}\cdots X_{j_r}^{\lambda_r}\]
where the sum is taken over all ordered $r$--tuples $(j_1,j_2,\ldots, j_r)$ with $1\le j_1<j_2<\cdots <j_r\le n$. We define $p^{\emptyset}:=1$, and we also write $p^{\lambda_1, \lambda_2, \ldots, \lambda_r}$ in lieu of $p^{\lambda}$ when $\lambda=(\lambda_1, \lambda_2, \ldots, \lambda_r)$.

For example, $p^{2,1}=X_1^2X_2+X_1^2X_3+X_2^2X_3$ and $p^{1,2}=X_1X_2^2+X_1X_3^2+X_2X_3^2$ so that $m_{2,1}=p^{2,1}+p^{1,2}$.

\subsection{The Hecke algebra}

The Hecke algebra $\H_n$ of $S_n$ is the unital associative algebra generated over $R=\Z[\xi]$ by the set $\{\T_s\mid s\in S\}$, {with relations}
\[
\T_w\T_s=\begin{cases}\T_{ws}+\xi\T_{w}&\text{if }\ell(ws)<\ell(w)\\ \T_{ws}&\text{otherwise,}\end{cases}
\]
where $\T_w:=\T_{s_{i_1}}\dots\T_{s_{i_r}}$ if $w=s_{i_1}\dots s_{i_r}$ is reduced. The Hecke algebra $\H_n$ is an $R$-module with basis $\{\T_w\mid w\in S_n\text{ reduced}\}$, and specializes at $\xi=0$ to the symmetric group algebra $\Z S_n$. We let $Z(\H_n)$, respectively $Z(\Z S_n)$, denote the centre of $\H_n$, respectively $\Z S_n$. Specialization of $Z(\H_n)$ at $\xi=0$ gives $Z(\Z S_n)$.
It is known that $Z(\H_n)$ has an integral $R$--basis $\{\Gamma_\lambda\mid \lambda\vdash n\}$ characterized by the following properties:
\begin{itemize}
\item $\Gamma_\lambda$ specializes at $\xi=0$ to the sum of elements in the conjugacy class $C_\lambda$ in $Z(\Z S_n)$, and
\item the only shortest elements from any conjugacy class appearing in $\Gamma_\lambda$ are those from $C_\lambda$, and they appear with coefficient 1.
\end{itemize}
These results can be found in \cite{GR97,Fmb}.
We refer to this basis in the sequel as the Geck--Rouquier basis.
If $\lambda=(\lambda_1,\lambda_2,\ldots ,\lambda_t)$, then we write $\Gamma_{\lambda}$ as $\Gamma_{\lambda_1,\lambda_2,\ldots ,\lambda_t}$.

The Hecke algebra is equipped with an inner product defined via the standard trace function $\tr (\sum_{w\in S_n}r_w\T_w)=r_1$ as follows,
\[\l \T_u,\T_v\r:=\tr(\T_v\T_u)=\begin{cases}1&\text{if }vu=1\\ 0&\text{otherwise.}\end{cases}\]
If $h$ is central in $\H_n$ then $\l \T_w,h\r$ is the coefficient of $\T_w$ in $h$.
This inner product gives rise, via specialization at $\xi=0$, to an inner product on the group algebra.

The Jucys--Murphy elements~\cite{Jucys1974,Mur83} $L_i$ are defined by setting $L_1=0$ and for $1\le i\le n-1$,
\[L_i:=\sum_{1\le j\le i-1}\T_{(j\ i)}.\]
For example, $L_2=\T_{s_1}$ and $L_3=\T_{s_1s_2s_1}+\T_{s_2}$. Specialization at $\xi=0$ gives corresponding Jucys--Murphy elements in $\Z S_n$.

Let $M_\xi=M_\xi(n)=\left\{m_\lambda(L_1,\dots ,L_n)\mid\ell(\lambda)\le n\right\}$, the set of monomial symmetric polynomials in  Jucys--Murphy elements in $Z(\H_n)$, and let $M_0$ denote the corresponding set in $Z(\Z S_n)$.

\section{Preliminaries}\label{sec:prelims}

In this section we recall some of the machinery from~\cite{FG:DJconj2006} necessary for the first main result of this paper, Theorem~\ref{thm:monbasis}.

For $k$ a positive integer, define
\[\a(k):=\sum_{m=1}^k \binom{k+m-1}{2m-1} \xi^{2m}.\]
This polynomial was introduced in~\cite{Mat99}, but can also be defined using $q$-integers, as in~\cite{FG:DJconj2006} (noting the different definition of $\xi$ used there).  For a composition $\lambda=(\lambda_1,\dots,\lambda_r)$, define $\a(\lambda)=\prod_{i=1}^r\a(\lambda_i)$.  For compositions $\lambda$ and $\mu$ {with $|\lambda|, |\mu|<n$}, define the polynomial $A_{\lambda,\mu}\in R$ to be the coefficient of $p^\lambda$ in $\sum_{|\gamma|=|\lambda|-|\mu|}\a(\gamma)p^\gamma p^\mu$, setting $A_{\lambda,\mu}=0$ if $|\lambda|\le|\mu|$.  Define the matrix $A^{(k)}$ by labelling its rows with compositions $\lambda$, where $|\lambda|<k$, according to the order defined in Section~\ref{Comp}, and with $(\lambda,\mu)$ entry {equal to} $A_{\lambda,\mu}$.  Define the matrix $Z^{(k)}$ recursively by setting $Z^{(1)}=(1)$ and
\[Z^{(k+1)}:=\begin{pmatrix} Z^{(k)}&I\\ 0&I\end{pmatrix}.\]
Define
\[\Xi^{(k+1)}:=\begin{pmatrix} \Xi^{(k)}&0\\ 0&\Xi^{(k)}\end{pmatrix}Z^{(k+1)}A^{(k+1)}\]
and $\Xi^{(1)}=(1)$.
The matrix $\Xi^{(k)}$ is the matrix of coefficients of minimal length conjugacy class elements in quasi-symmetric polynomials of Jucys--Murphy elements~\cite[Theorem 6.7]{FG:DJconj2006}.  Specifically,
\[
\langle\T_w,p^\mu(L_1,\dots,L_n)\rangle=\Xi^{(k)}_{\lambda',\mu'},
\]
for $\lambda,\mu\vDash k$ and $w$ of shape $\lambda$. For the sake of convenience, we use the bijection $\lambda \mapsto \lambda^{\prime}$, that maps compositions of a fixed positive integer $k$ to compositions of nonnegative integers less than $k$, to define the re-indexed matrix $X^{(k)}_{\lambda,\mu}:=\Xi^{(k)}_{\lambda',\mu'}$.  To obtain coefficients in symmetric polynomials (as opposed to quasi-symmetric polynomials) of Jucys--Murphy elements, define the matrix $T^{(k)}$ indexed by compositions $\lambda, \mu$ of $k$ by setting
  \[
  T_{\lambda,\mu}^{(k)}=\begin{cases}1&\text{if }\lambda=\mu\\ 1&\text{if }\hat\lambda=\mu\\ 0&\text{otherwise.}\end{cases}
  \]
Finally, $M^{(k)}$ is the matrix indexed by partitions of $k$ obtained from $(T^{(k)})^{-1}X^{(k)}T^{(k)}$ by deleting those rows and columns labelled by non-partitions.

For positive integers $n$ and $k\le n/2$, the $(\lambda,\mu)$--entry of $M^{(k)}$
is $\langle\T_w,m_\mu\rangle$, where $w\in S_n$ is an increasing element of shape $\lambda$ ~\cite[Theorem 7.1]{FG:DJconj2006} (there is always such a $w$ since $k\le n/2$). Also, by the nature of the construction of $M^{(k)}$, its entries are independent of $n$~\cite[Lemma 5.2]{FG:DJconj2006}.
One of the main results of~\cite{FG:DJconj2006} is that $M^{(k)}$ is  invertible over $R$.  This fact was originally conjectured by G.~James ~\cite{Mat99}. We give a brief summary of the construction of the inverse of $M^{(k)}$, needed in Section \ref{sec:basis}.

Let
\[\Upsilon^{(k+1)}:=Z^{(k+1)}A^{(k+1)}\begin{pmatrix}\Upsilon &0\\ 0&\Upsilon\end{pmatrix} \text{ with }\Upsilon^{(1)}=(1).\]
Then $\Xi$ and $K\Upsilon K$ are inverse ~\cite[Lemma 5.5]{FG:DJconj2006}, where $K$ is defined recursively by $K=\begin{pmatrix}K&-K\\ 0&-K\end{pmatrix}$. As with $\Xi$ above, we define the re-indexed matrix $Y_{\lambda,\mu}:=\Upsilon_{\lambda',\mu'}$.

The matrices $T^{-1}X T$ and $T^{-1}KYKT$ are inverse.
Let $N^{(k)}$ be the matrix obtained from $T^{-1}KYKT$ by deleting rows and columns not labelled by partitions.  It is shown in ~\cite{FG:DJconj2006} that $M^{(k)}N^{(k)}=I$, and the entries of $N^{(k)}$ are from $R$.

\section{A new integral basis for $Z(\H_n)$ using symmetric polynomials in Jucys--Murphy elements}\label{sec:basis}

In the case of the group algebra of the symmetric group, Murphy has shown that the set of monomial symmetric polynomials $m_\mu$ in Jucys--Murphy elements, such that an increasing element of shape $\mu$ is in $S_n$, is a $\Z$--basis for $Z(\Z S_n)$~\cite{Mur83}.  The corresponding statement for $Z(\H_n)$ appears in~\cite[Theorem 3.6]{Mat99}.
However, this generalization from $Z(\Z S_n)$ to $Z(\H_n)$ is incorrect.
We have verified that this particular set is not an $R$--basis for $Z(\H_n)$ when $n=3,4,5$, and we conjecture that this set is an $R$--basis for $Z(\H_n)$ only when $n=1,2$. We present the details here of the counterexample when $n=3$. In this case, the direct generalization of Murphy's symmetric group result states that the set of monomials $\left\{m_{\emptyset}, m_1,m_2\right\}$ in Jucys--Murphy elements is an $R$--basis for $Z(\H_3)$. If this is true, then the transition matrix from this monomial ``basis" to the Geck--Rouquier basis $\{\Gamma_{1,1,1},\Gamma_{2,1},\Gamma_{3}\}$ must be invertible in $R$. This transition matrix is:
\[\left[\begin{array}{ccc}
1&0&3\\
0&1&2\xi\\
0&0&1+\xi^2
\end{array}\right].\]
However, its determinant is 1+$\xi^2$, which is not invertible in $R$.
Observe that this set of monomials does specialize at $\xi=0$ to a $\Z$--basis for $Z(\Z S_3)$.

While a more general statement --- that there is an integral basis for $Z(\H_n)$ consisting solely of elements from $M_\xi$ --- also seems unlikely in general (though true for $n\le 4$, see Section~\ref{sec:monbases}), we show here that there always exists an integral basis for $Z(\H_n)$ which consists of $R$--linear combinations of elements from $M_\xi$. This is the first such basis for $Z(\H_n)$ given in terms of symmetric polynomials, and the second description of an integral basis after~\cite{GR97,Fmb}.

For a partition $\lambda$, let
\[
\M_\lambda=\sum_{\mu\vdash |\lambda|} N^{(|\lambda|)}_{\mu,\lambda}m_\mu,
\]
where $N^{(|\lambda|)}_{\mu,\lambda}$ is the entry in the matrix $N^{(|\lambda|)}$ whose row is labelled by  $\mu$ and whose column is labelled by $\lambda$.

\begin{thm}\label{thm:monbasis}
  The set $\B=\{\M_\lambda\mid |\lambda|+\ell(\lambda)\le n\}$ is an $R$--basis for $Z(\H_n)$.
\end{thm}

\begin{proof}
We will show that $\M_\lambda$ is the class element $\Gamma_{\overline\lambda}$ of the Geck--Rouquier basis, where $\overline\lambda$ is the unique partition of $n$ whose corresponding conjugacy class has a minimal length element of shape $\lambda$, plus a linear combination of class elements whose shortest elements are strictly shorter than $|\lambda|$.

Firstly, the coefficient of $\Gamma_\lambda$ in any central element $h$ is $\l \T_{w},h\r$ for $w$ an increasing element of shape $\lambda-1$, as a direct result of the characterization given in Section~\ref{sec:defs}.

For the moment, consider $n\ge 2k$ --- large enough so that all $w$ of shape $\lambda\vdash k$ are in $S_n$.  The elements of $M_\xi$, indexed by partitions ordered according to~\eqref{eq:def:compn.order}, are related to the class elements $\Gamma_\lambda$ via a block upper-triangular matrix, thanks to~\cite[Theorem 2.7]{Mat99}. That is, for monomial symmetric polynomials of degree less than or equal to $k$, we have
\[
(m_\emptyset,m_1,\dots,m_{1^k})=(\Gamma_{\overline\emptyset},\dots,\Gamma_{\overline{1^k}} )M,
\]
where $M$ is a block upper triangular matrix with diagonal blocks $M^{(0)}$, \dots, $M^{(k)}$. As noted in the previous section, the diagonal blocks are independent of $n$.

Restricting to monomials of degree $k$ we have
 \begin{multline*}
 (m_k,\dots,m_{1^k})=(\Gamma_{\overline{k}},\dots,\Gamma_{\overline{1^k}} )M^{(k)}%\\
 +(\textrm{linear combinations of $\Gamma_\nu$ for $|\nu-1|<k$.})
 \end{multline*}
Therefore
 \begin{multline*}
 (\M_{(k)},\dots,\M_{(1^k)})=(m_k,\dots,m_{1^k})N^{(k)}=(\Gamma_{\overline{k}},\dots,\Gamma_{\overline{1^k}} )\\
 +(\textrm{linear combinations of $\Gamma_\nu$ for $|\nu-1|<k$.})
 \end{multline*}
Because the entries of the matrices $M^{(k)}$ and $N^{(k)}$ are independent of $n$, we have
\[\M_\lambda=\Gamma_{\overline\lambda}+\textrm{a linear combination of $\{\Gamma_\nu\mid |\nu-1|<|\lambda|\}$.}
\]

The theorem follows by induction on $|\lambda|$.
\end{proof}

\begin{rem}
The bases given in Theorem \ref{thm:monbasis} embed in each other. That is, if $|\lambda|+\ell(\lambda)\le m\le n$, then the linear combination of elements from $M_\xi$ represented by $\M_\lambda$ is in the integral basis for both $Z(\H_m)$ and $Z(\H_n)$.

Theorem \ref{thm:monbasis} also provides an effective algorithm for computing an integral basis for $Z(\H_n)$.
The existing algorithm for computing class elements, given in~\cite{Fmb}, requires the computation of each class element independently, with terms $\T_w$ added individually.  That is, many computations are required for every conjugacy class in the symmetric group. However, the algorithm from Theorem~\ref{thm:monbasis} requires only the computation of the $n$ matrices $N^{(0)}, \ldots ,N^{(n-1)}$. Once these matrices are computed, basis elements corresponding to minimal elements of length less than $n$ can be obtained immediately in terms of elements from $M_\xi$. The problem of obtaining the matrices $N^{(k)}$ reduces to calculating the matrix $A^{(k)}$, defined in Section \ref{sec:prelims}, whose entries are determined by structure constants in the ring of quasi-symmetric polynomials. But most importantly, because of the recursive nature of the algorithm described in Theorem \ref{thm:monbasis}, once a basis of $Z(\H_{n-1})$ is obtained, only one additional matrix, $N^{(n-1)}$, is required to complete the basis for $Z(\H_n)$.
\end{rem}

To illustrate Theorem \ref{thm:monbasis}, we provide the bases for $Z(\H_n)$ produced therein when $n=3$, 4 and 5. The matrices $M^{(k)}$ and $N^{(k)}$, when $k\le 4$, are used in the calculations. For $k\le 3$ they are:
\[M^{(0)}=N^{(0)}=M^{(1)}=N^{(1)}=(1),\]
\[M^{(2)}=\begin{pmatrix}
  1+\xi^2&1\\ \xi^2&1
\end{pmatrix},\quad
N^{(2)}=\begin{pmatrix}
  1&-1\\ -\xi^2&1+\xi^2
\end{pmatrix},\]
\[
M^{(3)}=\begin{pmatrix}
  1+5\xi^2+5\xi^4+\xi^6 & 3+5\xi^2+\xi^4 & 1\\
  2\xi^2+4\xi^4+\xi^6   & 1+4\xi^2+\xi^4 & 1\\
  3\xi^4+\xi^6          & 3\xi^2+\xi^4   & 1
\end{pmatrix}
\]
and
\[
N^{(3)}=\begin{pmatrix}
  1+\xi^2       & -2\xi^2-3             & 2+\xi^2\\
  -2\xi^2-\xi^4 & 1+5\xi^2+2\xi^4       & -1-3\xi^2-\xi^4\\
  3\xi^4+\xi^6  & -3\xi^2-7\xi^4-2\xi^6 & 1+3\xi^2+4\xi^4+\xi^6
\end{pmatrix}.
\]
\begin{ex}\label{n=3}
$n=3$
\[\B=\{m_\emptyset,\ m_{1},\ m_{2}-\xi^2m_{1,1}\}.\]
Note that the last element in this basis is obtained from first column of $N^{(2)}$, corresponding to the partition $(2)$, because the element of shape $(2)$, $s_1s_2$, is in $S_3$.  The other column of $N^{(2)}$ is labelled by $(1,1)$, but $S_3$ contains no element of shape $(1,1)$.
\end{ex}
\begin{ex}\label{n=4}
$n=4$
\begin{align*}
  \B=\{&m_\emptyset,\ m_{1},\ m_2-\xi^2m_{1,1},\  -m_2+(1+\xi^2)m_{1,1},\ \\
&(1+\xi^2)m_{3}+(-2\xi^2-\xi^4)m_{2,1}+(3\xi^4+\xi^6)m_{1,1,1}\}.
\end{align*}

Here we have both partitions of $2$ giving an element in $S_4$, and so both columns of $N^{(2)}$ are represented.  We also have the first column of $N^{(3)}$ represented, labelled by the partition $(3)$.
\end{ex}
\begin{ex}\label{n=5}
$n=5$
\begin{align*}
  \B&=\{m_\emptyset,\ m_{1},\ m_2-\xi^2m_{1,1},\ -m_2+(1+\xi^2)m_{1,1},\\
&\quad (1+\xi^2)m_3-(2\xi^2+\xi^4)m_{2,1}+(3\xi^4+\xi^6)m_{1,1,1},\\
&\quad (-2\xi^2-3)m_3+(1+5\xi^2+2\xi^4)m_{2,1}+(-3\xi^2-7\xi^4-2\xi^6)m_{1,1,1},\\
&\quad (1+5\xi^2+5\xi^4+\xi^6)m_4-(\xi^8+6\xi^6+9\xi^4+4\xi^2)m_{2,2}-(\xi^{8}+6\xi^6+9\xi^4+3\xi^2)m_{3,1} \\
&\quad +(\xi^{10}+7\xi^{8}+14\xi^6+8\xi^4)m_{2,1,1}
-(\xi^{12}+8\xi^{10}+20\xi^{8}+16\xi^6)m_{1,1,1,1} \}.
\end{align*}
In this case, the first two columns of $N^{(3)}$, labelled by partitions $(3)$ and $(2,1)$, are represented, but not the third column, because $S_5$ contains no element of shape $(1,1,1)$. The final element of the basis is obtained from the first column of $N^{(4)}$, which corresponds to the partition $(4)$ and the increasing element $s_1s_2s_3s_4$ in $S_5$ ($N^{(4)}$ not shown).
\end{ex}

\section{Bases for $Z(\H_n)$ that are subsets of $M_\xi$}\label{sec:monbases}

While Theorem~\ref{thm:monbasis} gives a basis for $Z(\H_n)$ in terms of elements from $M_\xi$, it is clearly not a basis consisting solely of elements from $M_\xi$ when $n\ge 3$. The question of whether such a basis exists in general is open.

\begin{openquestion}\label{OQ1}
  For any $n$, does there exist an integral basis for $Z(\H_n)$ consisting solely of elements from $M_\xi$?
\end{openquestion}

When $n\le 4$, there is such a basis. For example, $Z(\H_3)$ has the basis $\{m_\emptyset,m_{1},m_{1,1}\}$, and $Z(\H_4)$ has bases \[\{m_\emptyset,m_{1},m_{2},m_{1,1},m_{1,1,1}\},  \{m_\emptyset,m_{1},m_{1,1},m_{1,1,1}, m_{2,1,1}\}\] \[\mbox{ and } \{m_\emptyset,m_{1},m_{1,1},m_{1,1,1}, m_{2,2,2}\}.\]

We do not know of any such basis of monomials for $Z(\H_5)$.
However, since we have a basis of monomials for $Z(\H_4)$, we would have an affirmative answer to Open Question \ref{OQ1} for $Z(\H_5)$ if we could find three monomials indexed by partitions of $n\ge 4$ that span the three class elements whose shortest terms are of length three or four ($s_1s_2s_3$, $s_1s_2s_4$ and $s_1s_2s_3s_4$), namely $\Gamma_{4,1}$, $\Gamma_{3,2}$ and $\Gamma_5$, modulo shorter class elements.  One of these could be $m_{1,1,1,1}$, which is equal to $\Gamma_5$.  Thus, the problem reduces to finding two monomials that span $\Gamma_{4,1}$ and $\Gamma_{3,2}$ modulo other class elements. Since there is no bound on the size of the partition for a monomial, beyond requiring the partition to have fewer than 5 parts (because the monomial is in $L_1,\dots,L_5$ and $L_1=0$), there are infinitely many candidates for such monomials.

\section{There is only one integral basis for $Z(\H_3)$ that is a subset of $M_\xi$} \label{sec:S3.monbasis}

Any subset of $M_\xi$ that is a basis for $Z(\H_n)$
specializes to a subset of $M_0$ that is a basis for $Z(\Z S_n)$.
Unfortunately, not every integral basis for $Z(\Z S_n)$ which consists of elements of $M_0$
is a specialization at $\xi=0$ of an integral basis for $Z(\H_n)$ which consists of elements of $M_\xi$.
See the counterexample in Section \ref{sec:basis}.  However, determining the
subsets of $M_0$ that are bases for $Z(\Z S_n)$ allows us to restrict our attention to the corresponding
subsets of $M_\xi$ in $Z(\H_n)$. In this section, we show that there are exactly four subsets of $M_0$ which are bases for $Z(\Z S_3)$, only one of which corresponds to an integral basis for $Z(\H_3)$.

\subsection{Coefficients of class sums in elements of $M_0$ for $\Z S_3$}

Let $\C$ be the matrix with three rows and infinitely many columns whose entries are the coefficients of the class elements $\Gamma_{\lambda}$ (rows) in the monomial symmetric polynomials $m_{\mu}$ (columns) in Jucys--Murphy elements for $\Z S_3$, where $\mu$ has fewer than three parts. In other words, the entries of $\C$ are $\l w_{\lambda-1}, m_{\mu} \r$, where $w_{\lambda-1}$ is an increasing element of shape $\lambda-1$, and $\lambda=(1,1,1)$, $(2,1)$, or $(3)$. For example, from Table \ref{tab:coeffs.S3}, which gives the first twenty columns of $\C$, the coefficient of $\Gamma_{2,1}$ in $m_{4,3}$ is $\l w_{1},m_{4,3}\r=\l s_1,m_{4,3}\r=8$. The rows of $\C$ are indexed by class elements ordered by the length of a shortest term in the class sum. The columns of $\C$ are indexed by monomial symmetric polynomials ordered from left to right according to the ordering on compositions given in Section \ref{Comp}. The monomials corresponding to partitions with more than two parts are not listed since in those situations the monomials are either zero (exactly three parts) or nonexistent (more than three parts).

For any integer $k$, we define the \emph{$k$--block} of $\C$ in the following way. If $k\ge 2$, then the $k$--block is the  submatrix of $\C$ consisting of those columns labelled by $m_\lambda$ for $\lambda\vdash k$. The $0$, or empty block, has the single column headed by $m_\emptyset$, while the $1$--block has the single column headed by $m_1$. A block is called \emph{even} or \emph{odd} according to whether $k$ is even or odd.

\begin{table}[h]
\caption[]{Coefficients of class sums in monomial symmetric polynomials of Jucys--Murphy elements in $\Z S_3$. This table was generated using GAP~\cite{Sch95} with CHEVIE~\cite{GHLMP}.}\label{tab:coeffs.S3}
\footnotesize
\begin{tabular}{lc@{\hspace{3mm}}c@{\hspace{3mm}}c@{\hspace{3mm}}c@{\hspace{3mm}}c@{\hspace{3mm}}c@{\hspace{3mm}}c@{\hspace{3mm}}c@{\hspace{3mm}}c@{\hspace{3mm}}c@{\hspace{3mm}}c}
&$m_\emptyset$&$m_{1}$&$m_{2}$&$m_{1,1}$&$m_{3}$&$m_{2,1}$&$m_{4}$&$m_{2,2}$&$m_{3,1}$\\
\hline
$\Gamma_{1,1,1}$       &1& 0& 3& 0& 0& 0& 7& 2& 2\\
$\Gamma_{2,1}$   &0& 1& 0& 0& 3& 2& 0& 0& 0\\
$\Gamma_{3}$&0& 0& 1& 1& 0& 0& 5& 1& 4\\
\\
&$m_{5}$&$m_{3,2}$&$m_{4,1}$&$m_{6}$&$m_{3,3}$&$m_{4,2}$&$m_{5,1}$&$m_{7}$&$m_{4,3}$&$m_{5,2}$&$m_{6,1}$\\
\hline
$\Gamma_{1,1,1}$       & 0& 0& 0& 23&2& 8& 10& 0& 0& 0& 0\\
$\Gamma_{2,1}$   & 11&4& 6&  0& 0& 0& 0&43&8&12&22 \\
$\Gamma_{3}$ & 0& 0& 0& 21&3& 6& 12& 0& 0& 0& 0
\end{tabular}
\normalsize
\end{table}

Certain rows of zeros appearing in the $k-$blocks ``alternate" according to whether $k$ is even or odd. More precisely, we have the following:

\begin{lem}\label{lem:chequerboard}
Let $k$ be a nonnegative integer, and let $\mu \vdash k$. Let $w_{\lambda-1}$ be an element of shortest length appearing in $\Gamma_{\lambda}$. If $\ell (w_{\lambda-1}) \not \equiv k \pmod{2}$, then $\l w_{\lambda-1},m_{\mu} \r=0$.
\end{lem}

\begin{proof}
  A column from an even block is a sum of products of Jucys--Murphy elements; the products being of even total degree.  Each group element $w$ in a Jucys--Murphy element is of odd length, because it is a reflection.  A product of two elements of odd length is even, of even length is even, and of mixed lengths is odd.  Thus, the length of a product of an even number of odd-length elements will be even. Consequently, the coefficient of any element of odd length ($s_1$ and $s_1s_2s_3$) in such a product will be zero.

  Similarly, a column from an odd block is a sum of products of an odd number of Jucys--Murphy elements.  Consequently, the coefficient of even length elements (the identity and $s_1s_2$) will be zero.
\end{proof}

\begin{lem}\label{lem:mon.relations.S3} The following relations hold {for} monomial symmetric polynomials in two variables:
\begin{enumerate}
 \item $m_i=m_2m_{i-2}-m_{2,2}m_{i-4}$ for $i\ge 5$\label{rel1}
 \item $m_{i,i}=m_{1,1}m_{i-1,i-1}$ for $i\ge 2$\label{rel2}
 \item $m_{i+j,i}=m_{i,i}m_{j}$\label{rel3} for $i,j\ge 1$.
  \end{enumerate}
\end{lem}
\begin{proof}
Straightforward.
\end{proof}

Because there are only two non-zero Jucys--Murphy elements in $\Z S_3$, namely the specializations of $L_2$ and $L_3$, Lemma \ref{lem:mon.relations.S3} gives relations among elements of $M_0$ for $\Z S_3$.

\begin{lem}\label{lem:mon.recurrences.S3}
  The following monomial recursions hold in $Z(\Z S_3)$:
  \begin{enumerate}
  \item $m_i=5m_{i-2}-4m_{i-4}$ for $i\ge 5$\label{rec1}
  \item $m_{i,i}=m_{i-1,i-1}+2m_{i-2,i-2}$  for $i\ge 3$\label{rec2}
  \item $m_{i+j,i}=5m_{i+j-2,i}-4m_{i+j-4,i}$  for {$i\ge 1$ and $j\ge 5$.}\label{rec3}
  \end{enumerate}
 \end{lem}

\begin{proof}
  To prove \ref{rec1}, we first verify the recurrence for $i=5,6,7,8$.
  For instance, from Table~\ref{tab:coeffs.S3}, we have
  \[m_5=\begin{pmatrix}0\\ 11\\ 0\end{pmatrix}=5\begin{pmatrix}0\\ 3\\ 0\end{pmatrix}-4\begin{pmatrix}0\\ 1\\ 0\end{pmatrix}=5m_{3}-4m_{1}\]
  and
  \[m_6=\begin{pmatrix}23\\ 0\\ 21\end{pmatrix}=5\begin{pmatrix}7\\ 0\\ 5\end{pmatrix}-4\begin{pmatrix}3\\ 0\\ 1\end{pmatrix}=5m_{4}-4m_{2},\]
where each monomial is written as a vector with respect to the basis of class sums for the centre.

If $i\ge 9$, we have by Lemma \ref{lem:mon.relations.S3} \ref{rel1} and induction,
\begin{align*}
  m_i&=m_2m_{i-2}-m_{2,2}m_{i-4}\\
     &=m_2\left(5m_{i-4}-4m_{i-6}\right)-m_{2,2}\left(5m_{i-6}-4m_{i-8}\right)\\
     &=5\left(m_2m_{i-4}-m_{2,2}m_{i-6}\right)-4\left(m_2m_{i-6}-m_{2,2}m_{i-8}\right)\\
     &=5m_{i-2}-4m_{i-4}.
\end{align*}
The proof of \ref{rec2} is similar. To prove \ref{rec3}, we use Lemma \ref{lem:mon.relations.S3} \ref{rel3} and recurrence \ref{rec1} of this lemma:
\begin{align*}
  m_{i+j,i}&=m_{i,i}m_j\\
     &=m_{i,i}\left(5m_{j-2}-4m_{j-4}\right)\\
     &=5m_{i,i}m_{j-2}-4m_{i,i}m_{j-4}\\
     &=5m_{i+j-2,i}-4m_{i+j-4,i}.
\end{align*}
\end{proof}

\begin{lem}\label{lem:mi}
  If $i\ge 1$, then
  \begin{align*}
    \l 1,m_{i}\r&=\frac{1}{6}\left(1+(-1)^i\right)\left(2^i+5\right)\\
    \l s_1,m_i\r&=\frac{1}{6}\left(1-(-1)^i\right)\left({2^i+1}\right)\\
    \l s_1s_2,m_{i}\r&=\frac{1}{6}\left(1+(-1)^i\right)\left(2^i-1\right).
  \end{align*}
\end{lem}

\begin{proof}
 The lemma is verified easily for $1\le i\le 5$ using Table~\ref{tab:coeffs.S3}.  For $i>5$, we use the recurrences of Lemma~\ref{lem:mon.recurrences.S3}, and proceed by induction.

 For instance, Table~\ref{tab:coeffs.S3} gives, respectively for $i=1,\dots,5$, the values $0,3,0,7,0$ of the inner product $\l 1,m_{i}\r$, and this agrees with the claim in each case.  If $i>5$, then $m_i=5m_{i-2}-4m_{i-4}$ by Lemma~\ref{lem:mon.recurrences.S3}, and so

 \begin{align*}
    \l 1,m_i\r&=5\l 1,m_{i-2}\r-4\l 1, m_{i-4}\r\\
        &=5\left[\frac{1}{6}\left(1+(-1)^{i-2}\right)\left(2^{i-2}+5\right)\right]-4\left[\frac{1}{6}\left(1+(-1)^{i-4}\right)\left(2^{i-4}+5\right)\right]\\
        &=\frac{1}{6}\left(1+(-1)^{i}\right)\left[5\left(2^{i-2}+5\right)-4\left(2^{i-4}+5\right)\right]\\
        &=\frac{1}{6}\left(1+(-1)^{i}\right)\left[\left(5.2^2-4\right)2^{i-4}+5\right]\\
        &=\frac{1}{6}\left(1+(-1)^{i}\right)\left(2^i+5\right).
 \end{align*}
 The other inductions are similar.
\end{proof}

\begin{lem}\label{lem:mii,mij}
  If $i,j\ge 1$, then
   \begin{align*}
    \l 1,m_{i,i}\r&=\frac{1}{3}\left(2^i+2(-1)^{i}\right)\\
    \l s_1,m_{i,i}\r&=0\\
    \l s_1s_2,m_{i,i}\r&=\frac{1}{3}\left(2^i-(-1)^{i}\right)\\
     \l 1,m_{i+j,i}\r&= \frac{1}{6}\left(1+(-1)^j\right)\left(2^{i+j}+2^i+4(-1)^i\right)\\
    \l s_1,m_{i+j,i}\r&= \frac{1}{6}\left(1-(-1)^j\right)\left(2^{i+j}+2^i\right)\\
    \l s_1s_2,m_{i+j,i}\r&= \frac{1}{6}\left(1+(-1)^j\right)\left(2^{i+j}+2^i-2(-1)^i\right).
  \end{align*}
\end{lem}

\begin{proof}
The first three equalities are easily verified using Table~\ref{tab:coeffs.S3} when $1\le i \le 2$. For $i\ge 3$, use Lemma \ref{lem:mon.recurrences.S3} \ref{rec2} and proceed by induction.

For the second three equalities, we first verify the eight cases: $1\le i\le 2$, $1\le j\le 4$, all but one of which is contained in Table~\ref{tab:coeffs.S3}. Then, for fixed $i\ge 1$, the closed forms follow by induction on $j$, using Lemma \ref{lem:mon.recurrences.S3} \ref{rec3}. Similarly, for fixed $j\ge 1$, we use Lemma \ref{lem:mon.relations.S3} \ref{rel3} and Lemma \ref{lem:mon.recurrences.S3} \ref{rec2} to write
\[m_{i+j,i}=m_{i,i}m_{j}=(m_{i-1,i-1}+2m_{i-2,i-2})m_{j}=m_{i-1+j,i-1}+2m_{i-2+j,i-2},\]
and proceed by induction on $i$.
\end{proof}

\begin{cor}\label{cor:beyond-ii}
  All coefficients of class elements in the monomials $m_{i+j,i}$, for {$i,j\ge 1$}, are even.
\end{cor}

\subsection{Monomial bases for $Z(\Z S_3)$ and $Z(\H_3)$}
\begin{lem}
  The only pair of non-trivial monomials that integrally span $\{\Gamma_{1,1,1},\Gamma_{3}\}$ is $\{m_2,m_{2,2}\}$.
\end{lem}
\begin{proof}
  To begin with, Lemma~\ref{lem:chequerboard} tells us that any two-element integral spanning set for $\{\Gamma_{1,1,1},\Gamma_{3}\}$ cannot contain monomials from odd blocks.  Secondly, we are restricted to monomials of form $m_{2j}$ and $m_{i,i}$ in light of Corollary~\ref{cor:beyond-ii} (for $i,j\ge 1$). Lemmas \ref{lem:mi} and \ref{lem:mii,mij} imply that the matrix of coefficients for the rows labelled by $\Gamma_{1,1,1}$ and $\Gamma_{3}$ for the two monomials $m_{2j}$ and $m_{i,i}$ is
  \[\begin{pmatrix}
%  &\\
\left( 1+(-1)^{2j}\right) \left( 2^{2j}+5\right) /6 & &\left(2^i+2(-1)^{i}\right)/3\\
%&\\
&\\
\left( 1+(-1)^{2j}\right) \left( 2^{2j}-1\right)/6 & &\left(2^i-(-1)^{i}\right)/3\\
%&
  \end{pmatrix},\]
with determinant \[d=\frac{\left(-1\right)^{i+1}\left(4^j+1\right)+2^{i+1}}{3}.\]
If $d=\pm 1$, then $\{m_{2j}, m_{i,i}\}$ is an integral spanning set for $\{\Gamma_{1,1,1},\Gamma_{3}\}$. We determine all values for $i$ and $j$ such that $d=\pm 1$.

If $i$ is odd, then $d>0$, and so we only need to find values of $i$ and $j$ for $d=1$. But, since $i,j\ge 1$, it follows that $d\equiv 3 \pmod{4}$, and so there are no solutions in this case.
If $i$ is even, then $d\equiv 1 \pmod{4}$.
Thus, $d=1$, and we must solve the exponential Diophantine equation
\[2^{i+1}=4^j+4.\] If $j\ge 2$, then reduction modulo 8 shows that this equation has no solutions. Therefore, the only solution is $i=2$ and $j=1$, which gives the set of monomials $\{m_2,m_{2,2}\}$.
\end{proof}

Allowing the trivial monomial $m_\emptyset$ produces only three more pairs of monomials that integrally span $\{\Gamma_{1,1,1},\Gamma_{3}\}$:

\begin{prop}\label{prop:mons.spanning.1,s12}
  The only pairs of monomials that integrally span $\{\Gamma_{1,1,1},\Gamma_{3}\}$ are
 $\{m_\emptyset,m_2\}$, $\{m_\emptyset,m_{1,1}\}$, $\{m_\emptyset,m_{2,2}\}$ and $\{m_2,m_{2,2}\}$.
\end{prop}

\begin{proof}
  If the monomial $m_\emptyset$ is included in a spanning set, the other monomial included to span $\{\Gamma_{1,1,1},\Gamma_{3}\}$ must have coefficient one on $\Gamma_{3}$. The only monomials satisfying this condition are $m_2$, $m_{1,1}$ and $m_{2,2}$ (Lemma~\ref{lem:mi} and Lemma~\ref{lem:mii,mij}).
\end{proof}

\begin{prop}\label{prop:mon.spanning.s1}
  The only monomial that integrally spans $\{\Gamma_{2,1}\}$ is $m_1$.
\end{prop}
\begin{proof}
  To span this one-element set integrally, a monomial must clearly have coefficient 1 or $-1$ on $\Gamma_{2,1}$. From Lemma~\ref{lem:mi}, the only monomial of form $m_i$ that has coefficient 1 or $-1$ on $\Gamma_{2,1}$ is $m_1$.  Monomials from even blocks are excluded because they have coefficient zero on $\Gamma_{2,1}$ (Lemma~\ref{lem:chequerboard}).  At the same time, other monomials in odd blocks are all of form $m_{i,j}$ for $i\neq j$, and so always have even coefficients by Corollary~\ref{cor:beyond-ii}.  This proves the claim.
\end{proof}

\begin{thm}\label{thm:mon.bases.ZS3}
 There are exactly four bases for $Z(\Z S_3)$ consisting solely of {elements from $M_0$}: $\{m_\emptyset,m_1,m_2\}$, $\{m_\emptyset,m_1,m_{1,1}\}$, $\{m_\emptyset,m_1,m_{2,2}\}$ and $\{m_1,m_2,m_{2,2}\}$.
\end{thm}
\begin{proof}
  This follows immediately from Propositions~\ref{prop:mons.spanning.1,s12} and~\ref{prop:mon.spanning.s1}.
\end{proof}

\begin{thm}
  The only {subset of $M_\xi$ that is} an integral basis for $Z(\H_3)$ is $\{m_\emptyset,m_1,m_{1,1}\}$.
\end{thm}
\begin{proof}
  The proof reduces to checking the columns corresponding to the monomials that form bases for $Z(\Z S_3)$ (Theorem~\ref{thm:mon.bases.ZS3}). These are:
  \[
  \begin{array}{lccccc}
    &m_\emptyset&m_1&m_2&m_{1,1}&m_{2,2}\\
    \hline
    \Gamma_{1,1,1}&1&0&3&0&2+\xi^2\\
    \Gamma_{2,1}&0&1&2\xi&0&\xi(3+\xi^2)\\
    \Gamma_{3}&0&0&1+\xi^2&1&1+4\xi+\xi^2
  \end{array}
  \]
  The only three columns that together give an integrally invertible matrix are $m_\emptyset$, $m_1$ and $m_{1,1}$.
\end{proof}


\begin{thebibliography}{1}

\bibitem{DJ87}
Richard Dipper and Gordon James.
\newblock Blocks and idempotents of {H}ecke algebras of general linear groups.
\newblock {\em Proc. London Math. Soc. (3)}, 54(1):57--82, 1987.

\bibitem{Fmb}
Andrew Francis.
\newblock The minimal basis for the centre of an {I}wahori-{H}ecke algebra.
\newblock {\em J. Algebra}, 221(1):1--28, 1999.

\bibitem{FG:DJconj2006}
Andrew~R. Francis and John~J. Graham.
\newblock {Centres of Hecke algebras: the Dipper-James conjecture}.
\newblock {\em J. Algebra}, 306:244--267, 2006.

\bibitem{GHLMP}
M.~Geck, G.~Hiss, F.~L{\"u}beck, G.~Malle, and G.~Pfeiffer.
\newblock {CHEVIE}---a system for computing and processing generic character
  tables. {Computational} methods in {Lie} theory ({Essen} 1994).
\newblock {\em Appl.~Algebra ENGRG.~Comm.~Comput.}, 7(3):175--210, 1996.

\bibitem{GR97}
Meinolf Geck and Rapha{\"e}l Rouquier.
\newblock Centers and simple modules for {I}wahori-{H}ecke algebras.
\newblock In {\em Finite reductive groups (Luminy, 1994)}, pages 251--272.
  Birkh\"auser Boston, Boston, MA, 1997.

\bibitem{Jucys1974}
A.~A.~A. Jucys.
\newblock Symmetric polynomials and the center of the symmetric group ring.
\newblock {\em Rep. Mathematical Phys.}, 5(1):107--112, 1974.

\bibitem{Lascoux2006}
Alain Lascoux.
\newblock The {Hecke} algebra and structure constants of the ring of symmetric
  polynomials.
\newblock arXiv:math.CO/0602379 v1 17 Feb 2006.

\bibitem{Mat99}
Andrew Mathas.
\newblock Murphy operators and the centre of {Iwahori}-{Hecke} algebras of type
  {$A$}.
\newblock {\em J. Algebraic Combinatorics}, 9:295--313, 1999.

\bibitem{Mur83}
G.~E. Murphy.
\newblock The idempotents of the symmetric group and {N}akayama's conjecture.
\newblock {\em J.Algebra}, 81:258--265, 1983.

\bibitem{Sch95}
Martin Sch{\"o}nert et~al.
\newblock {\em {GAP} --- {Groups}, {Algorithms}, and {Programming}}.
\newblock Lehrstuhl D f{\"u}r Mathematik, Rheinisch Westf{\"a}lische Technische Hochschule, Aachen, Germany, fifth edition, 1995.

\end{thebibliography}
\end{document}